\documentclass[12pt,a4paper]{article}

\usepackage{amssymb}

\newtheorem{Proposition}{Proposition}
\newtheorem{Theorem}[Proposition]{Theorem}

\newtheorem{Lemma}[Proposition]{Lemma}
\newtheorem{Remark}[Proposition]{Remark}
\newtheorem{Example}[Proposition]{Example}

\newcommand{\C}{{\mathbb C}}
\newcommand{\N}{\mathbb N}
\newcommand{\Z}{\mathbb Z}

\newcommand{\Qp}{\mathbb Q_p} 
\newcommand{\Qps}{\mathbb Q_{p^s}} 

\newcommand{\Q}{\mathbb Q} 
\newcommand {\F}{\mbox{${\cal F}$}}
\newcommand {\Fl}{\mbox{${\cal F} \!\!\!\!\!\!\!\; {\cal J}$}}
\newcommand{\proof}{\noindent {\bf Proof: }}
\newcommand {\qed}{\hfill $\square$}

\renewcommand{\P}{\mathbb P} 

\begin{document}

\title{On Harder-Narasimhan strata in flag manifolds}
\author{Sascha Orlik}
\date{}
\maketitle

\begin{abstract}
This paper deals with a question of Fontaine and Rapoport which was posed in \cite{FR}.
There they asked for the determination of the index set of the Harder-Narasimhan vectors  of the filtered isocrystals
with fixed  Newton- and  Hodge vector.
The aim of this paper is to give an answer to their question.
\end{abstract}

\section{Introduction}
Let $k$ be an algebraically closed  field of characteristic $p>0$ and
denote by $K$=Quot($W(k)$) the fraction field of the ring
of Witt vectors. Let $\sigma \in  \mbox{Aut}(K/\mathbb Q_p)$ be the Frobenius automorphism. 
Fix an isocrystal $(V,\Phi)$ of dimension $d$ over $k,$ i.e., $V$ is a $d$-dimensional vector space 
over $K$ together with a $\sigma$-linear bijective endomorphism $\Phi$ of $V.$
Consider a  filtration $\F$ of $V$ by subspaces, which yields the notion of a filtered isocrystal $(V,\Phi,\F)$ over $k$. 
Analogously to the case of vector bundles on curves one can  define   numerical invariants as  the degree
and the slope for these objects. These numerical invariants lead by the usual machinery to the definition of a  
semistable filtered isocrystal and 
to the Harder-Narasimhan filtration of a filtered isocrystal.
To every HN-filtration one  associates its HN-vector 
$$\lambda=(\lambda_1,\dots,\lambda_d)\in (\Q^d)_+:=\{(v_1,v_2,\ldots, v_d)\in \Q^d; v_1\geq v_2 \geq \ldots \geq v_d \}.$$
This vector describes the HN-polygon attached to the HN-filtration of the object $(V,\Phi,\F).$ 
A natural question, raised in \cite{FR}, is to determine the set of HN-vectors 
for a fixed isocrystal $(V,\Phi)$ and  a fixed Hodge vector $\mu=(\mu_1,\ldots,\mu_d) \in (\Q^d)_+$
describing the type of the considered filtrations $\F.$
The present paper gives an answer to this problem. 

Let $\nu=\nu(V,\Phi)\in (\Q^d)_+$ be the Newton vector of the isocrystal $(V,\Phi).$ Then
we may write $\nu$ as
$\nu=(\nu(1)^{s_1},\nu(2)^{s_2},\ldots,\nu(l)^{s_l})$
(the exponents indicate that the term $\nu(i)$ is repeated exactly $s_i$ times),  
with
$\nu(1) \geq \nu(2) \geq \dots \geq \nu(l), \mbox{ where } \nu(i)=\frac{r_i}{s_i}\in \Q$
for some relatively prime integers $r_i,s_i \in \Z, s_i > 0.$ As in \cite{FR} we put
$|v|:=\sum_{i=1}^d v_i$ for any $v\in (\Q^d)_+.$ Further we consider on $(\Q^d)_+$ an
order $\leq$ which is a slight generalisation of the usual dominance order on $(\Q^d)_+$ in loc.cit.
Viewing tuples in $(\Q^d)_+$ as cocharacters of $GL_d,$ the order  $\leq$ is induced by
the projection $GL_d \rightarrow PGL_d$ and the dominance order on $PGL_d$ (compare the end of section 2).
\begin{Theorem}
Let $\lambda=(\bar{\lambda}_1^{d_1},\dots,\bar{\lambda}_r^{d_r})\in (\Q^d)_+$ with
$\bar{\lambda}_1 > \dots > \bar{\lambda}_r$ and $d_i\geq 1$ for $i=1,\ldots,r.$ Then $\lambda$ appears as the HN-vector 
of a filtration $\F$ on $V$ of type $\mu$ if and only if 
there exist disjoint partitions into non-empty subsets of the intervals $[1,l]=I_1\cup \ldots \cup I_r$ and
$[1,d]=J_1\cup \ldots \cup J_r,$ such that the following conditions are satisfied for all $i=1,\ldots,r.$
\begin{enumerate}
\item $|J_i|=d_i,$ 
\item $d_i=\sum_{k\in I_i}s_k ,$

If we put $\underline{\nu}_i=((\frac{r_k}{s_k})^{s_k}; k\in I_i)_+ \in (\Q^{d_i})_+$ and 
$\underline{\mu}_i=(\mu_k; k\in J_i)_+ \in (\Q^{d_i})_+$ (the notation $\underline{\mu}_i=(\mu_k; k\in J_i)_+$ means that the entries $\mu_k$
are ordered with multiplicities in a decreasing tuple) then 
\item $\underline{\mu}_i \geq \underline{\nu}_i,$
\item $ \bar{\lambda}_i\cdot d_i =|\underline{\mu}_i|-|\underline{\nu}_i|.$
\end{enumerate}
\end{Theorem}
The strategy for proving this theorem works as follows.
Let $K$=Quot($W(k)$) be the fraction field of the ring of Witt vectors. We will
construct a stratification of the adic space $\Fl^{ad}$ (\cite{H}) associated to the flag variety $\Fl=\Fl(V,\mu)$ 
over $K$ consisting of flags of type $\mu.$ We will prove:
\begin{Theorem}
There exists a stratification $\Fl^{ad}=\bigcup\limits^\cdot_{\gamma \in \Gamma}\Fl^{ad}_\gamma$ by locally closed pseudo-adic subspaces.
Each stratum is a vector bundle over a continuous family of products of certain period domains (confer to chapter 3 for the notion of a period domain).
\end{Theorem}
A stratum is characterized by the structure of the Harder-Narasimhan filtration attached to the objects $(V,\Phi,\F).$
Therefore the stratification is finer than that one induced by the HN-vectors.
In the basic isocrystal case it coincides with the stratification considered by Kottwitz and Rapoport on their work on the Euler-Poincar\'e characteristic
of period domains (compare \cite{R2}, \cite{R4}). Since we have an explicit description of the index set $\Gamma$ we are able to
prove Theorem 1.

The structure of this paper is given as follows. In section 2 we introduce some notations. Further we determine those Newton- and Hodge vectors
for which there exist semistable filtrations. Section 3 deals with the stratification mentioned in Theorem 2. 
In section 4 we turn to the proof of Theorem 1.

I would like to thank M. Rapoport for telling me about the problem on the HN-vectors. I am especially grateful for his steady interest and his remarks 
on this paper. I also would like to thank A. Huber and  B. Totaro for helpful comments on this paper.

\section{The existence of semistable filtrations}
Let $k$ be an algebraically closed  field of characteristic $p>0.$ We
denote by $K$=Quot($W(k)$) the fraction field of the corresponding  ring
of Witt vectors $W(k).$
Fix an isocrystal $(V,\Phi)$ of dimension $d$ over $k.$ 
Denote by $J=Aut(V,\Phi)$ its automorphism group (compare \cite{RR}). It is an algebraic reductive group
defined over $\Qp.$
Let 
$$\nu=\nu(V,\Phi) \in (\Q^d)_+:=\{(v_1,v_2,\ldots, v_d)\in \Q^d; v_1\geq v_2 \geq \ldots \geq v_d \}$$
be the Newton vector of $(V,\Phi)$ [loc.cit]. Then $\nu$ has the form
$$\nu=(\underbrace{\nu(1),\ldots,\nu(1)}_{s_1},\underbrace{\nu(2),\ldots,\nu(2)}_{s_2},\ldots,\underbrace{\nu(l),\ldots,\nu(l)}_{s_l})$$ 
with
$$\nu(1) \geq \nu(2) \geq \dots \geq \nu(l), \mbox{ where } \nu(i)=\frac{r_i}{s_i}\in \Q$$
for some relatively prime integers $r_i,s_i \in \Z, s_i > 0.$
The tuples $$(\underbrace{\nu(i),\ldots,\nu(i)}_{s_i}),\; i=1,\ldots,l$$ are just the
Newton vectors of the simple summands of $(V,\Phi).$ We fix a decomposition of $(V,\Phi)$ into simple
subisocrystals $V=\oplus_{i=1}^l V_i$ with $slope(V_i)=\nu(i).$ We denote by $V^\bullet$ 
the slope filtration of $(V,\Phi).$ It is defined by
$$V^x=\oplus_{i\leq -x} V_i.$$
Of course this definition does not depend on the chosen decomposition.
We stress that none of the results of this paper depends on our chosen decomposition, since any other 
one is obtained by applying an element of $J(\Qp)$ to the fixed one.

Fix a Hodge vector
$$\mu=(\mu_1 \geq \mu_2 \geq \ldots \geq \mu_d) \in (\Q^d)_+ .$$ 
Denote by $\Fl=\Fl(V,\mu)/K$ the flag variety
consisting of flags of type $\mu.$ This variety is defined over $\Qp$. If we rewrite $\mu$ as
$$\mu=(\bar{\mu}_1^{g(\bar{\mu}_1)},\ldots,\bar{\mu}_s^{g(\bar{\mu}_s)}),$$
where $\bar{\mu}_1 > \cdots >\bar{\mu}_s$ then 
$$\Fl(K)=\{\Q\mbox{-filtrations $\F$  on $V$} ; \dim gr^{\bar{\mu}_i}_{\F}(V) = g(\bar{\mu}_i), \;i=1,\ldots,s \}.$$
Inside $\Fl(K)$ we have the - possibly empty - subset $\Fl^{wa}=\Fl(V,\Phi,\mu)^{wa}$ of weakly admissible filtrations. 
A filtration $\F \in \Fl(K)$ is called weakly admissible if the filtered isocrystal $(V,\Phi,\F)$ is weakly admissible, i.e., if
$$slope_{\F}(U) \leq slope_{\F}(V)$$ 
for all subisocrystals  $U\subset V$  (i.e., $\F$ is semistable ) and
$$slope_{\F}(V)=0.$$
Here we put 
$$slope_{\F}(U):=(\deg_{\F}(U) + \deg_{V^\bullet}(U))/ \dim V,$$ 
for any subisocrystal $U\subset V,$
where $\deg_{\F}U = \sum_{x \in \Q} x \dim gr^x_{\F}(U)$ and $\deg_{V^\bullet}U = \sum_{x \in \Q} x \dim gr^x_{V^\bullet}(U).$
By \cite{FR} we know that $\Fl^{wa} \neq \emptyset$ if and only if the following holds:
\begin{eqnarray}
\mu_1 & \geq & \nu_1 \nonumber\\
\mu_1 + \mu_2 & \geq & \nu_1 + \nu_2 \nonumber\\
&\vdots& \\
\mu_1 + \ldots + \mu_{d-1} &\geq& \nu_1 + \ldots +  \nu_{d-1} \nonumber\\
\mu_1 + \ldots + \mu_d &=& \nu_1 + \ldots +  \nu_d \nonumber
\end{eqnarray}
As in \cite{FR} we abbreviate this system of inequalities in writing $\mu \geq \nu.$ Further we put 
$|v|:=v_1 + \ldots + v_d \in \Q$ for any $v \in \Q^d.$
We denote for every $\F \in \Fl(K)$  by $HN(\F)$ the Harder-Narasimhan filtration (\cite{RZ} Prop. 1.4)
$$(0)=W^0 \subset W^1 \subset W^2 \subset \cdots \subset W^r=V$$ 
of the filtered isocrystal $(V,\Phi,\F).$ This filtration is uniquely determined by the property that
the induced filtered isocrystal $gr^i(HN(\F))$ is semistable for all $i$ and that
$$slope_{\F}(gr^1(HN(\F))) > slope_{\F}(gr^2(HN(\F))) > \cdots > slope_{\F}(gr^r(HN(\F))).$$
Put 
$$|HN(\F)|:=r.$$
Thus an element $\F\in \Fl(K)$ lies in $\Fl^{wa}$ if and only if $|HN(\F)|=1$ and $slope_{\F}(V)=0.$
We denote by $\Fl^{ss}=\Fl(V,\Phi,\mu)^{ss} \subset \Fl(K)$ the subset of semistable filtrations. Then we can prove:

\begin{Lemma}
The subset $\Fl^{ss}$ is non-empty if and only if  
\begin{eqnarray}
\mu_1 + \frac{|\nu|}{d} & \geq & \nu_1 + \frac{|\mu|}{d} \nonumber\\
\mu_1 + \mu_2 + \frac{2|\nu|}{d} & \geq & \nu_1 + \nu_2 + \frac{2|\mu|}{d} \\
&\vdots & \nonumber\\
\mu_1 + \ldots +\mu_{d-1} + \frac{(d-1)|\nu|}{d} & \geq & \nu_1 + \ldots + \nu_{d-1} + \frac{(d-1)|\mu|}{d} \nonumber
\end{eqnarray}
\end{Lemma}
\proof Substituting the Hodge vector $\mu$ by  $\mu(\alpha):=\mu + \alpha(1,\ldots,1) \in (\Q^d)_+, \; \alpha \in \Q,$ does not change 
the set of  semistable  filtrations.
Indeed, if we denote by $slope^\alpha$ the modified slope function with respect to $\mu(\alpha)$ then we have 
$$ slope_{\F}^\alpha(U) = slope_{\F}(U)+\alpha$$
for any subisocrystal $U$ of $V.$  Put $\alpha:=\frac{|\nu| - |\mu|}{d} \in \Q.$ 
Then $|\mu(\alpha)| = |\nu|$ and thus $\Fl(V,\Phi,\mu(\alpha))^{wa}=\Fl(V,\Phi,\mu(\alpha))^{ss}=\Fl(V,\Phi,\mu)^{ss}.$ Using (1)
we conclude $\Fl^{ss} \neq \emptyset \Leftrightarrow$
\begin{eqnarray*}
\mu(\alpha)_1 & \geq & \nu_1 \\
\mu(\alpha)_1 + \mu(\alpha)_2 & \geq & \nu_1 + \nu_2 \\
&\vdots& \\
\mu(\alpha)_1 + \ldots + \mu(\alpha)_{d-1} &\geq & \nu_1 + \ldots +  \nu_{d-1} \\
|\mu(\alpha)| &=& |\nu| 
\end{eqnarray*}
Since $\mu(\alpha)_1+ \cdots + \mu(\alpha)_i= \mu_1+ \cdots +\mu_i + i \alpha$ the lemma is proved. \qed

\noindent As the inequalities of the lemma match with the relation $\mu \geq \nu$ in the case  where $|\mu| =|\nu|$, we may extend the order
to all of $\Q^d.$ Hence, we also write $\mu \geq \nu$ for two arbitrary $\mu,\nu \in \Q^d$ if the 
system of inequalities (2) is satisfied. 
Let $\mu \in \Q^d$ and view it as a cocharacter of $GL_d.$ If we want to express its image under the projection $GL_d \rightarrow PGL_d$
as a linear combination of the standard simple roots, then we just have to substitute $\mu$ by the tuple $\mu-\frac{|\mu|}{d}(1,\ldots,1) \in \Q^d.$
Thus we see that this new order corresponds to the $PGL_d$-dominance order on $\Q^d.$

\section{The stratification}
Let $$W^\bullet=((0)=W^0 \subset W^1 \subset W^2 \subset \cdots \subset W^r=V)$$ be a filtration of $V$
by subisocrystals. Our next goal is to examine when such a filtration appears as the HN-filtration of a flag $\F\in\Fl(K).$ 
Again, after applying a transformation in $J(\Qp),$ we may suppose without loss of generality that $W^\bullet$ consists 
of the pieces $V_i$ of our decomposition $V=\oplus_{i} V_i.$ 
Put for $1\leq i \leq r$
$$d_i:=\dim gr^i(W^\bullet), $$ 
and set $d_0=0.$ We denote for each $i$ by $\nu(W^\bullet,i):=\nu(gr^i(W^\bullet))\in (\Q^{d_i})_+$ the Newton vector
of the isocrystal $gr^i(W^\bullet).$ The entries of these vectors are induced by the 
Newton vector $\nu.$ In the following we want to make  precise this relationship .

We denote for an integer $n\in \N$ by $S_n$  the symmetric group of the set $\{1,\ldots,n\}.$ Then we may view $S_l$ as a subgroup of $S_d$
by permuting the blocks 
$$\{1, \ldots, s_1\},\{s_1+1 , \ldots ,s_1+s_2\},\ldots ,\{s_1+\ldots +s_{l-1}+1,\ldots,d\}$$ 
and leaving the entries within a block invariant. 
Via this embedding the induced action of $S_l$ on $\Q^d$ is given as follows.
Let
$$\lambda=(\underline{\lambda}_1,\ldots,\underline{\lambda}_l) \in \Q^d$$ with
$\underline{\lambda}_i \in \Q^{s_i},\; i=1,\ldots,l,$ and let $x\in S_l.$ Then we have
$$x\cdot \lambda=(\underline{\lambda}_{x^{-1}(1)},\ldots,\underline{\lambda}_{x^{-1}(l)}) \in \Q^d.$$
Especially if $\lambda=\nu$ then 
$$x\cdot\nu=(\nu({x^{-1}(1)})^{s_{x^{-1}(1)}},\ldots,\nu({x^{-1}(l)})^{s_{x^{-1}(l)}}).$$ 
Here and in the following it should be clear from the context, whether we consider an element
$x\in S_l$ as a permutation in $S_l$ or $S_d.$
As in the case of our Hodge vector $\mu$ we rewrite $\nu$ as 
$$\nu=(\bar{\nu}_1^{h(\bar{\nu}_1)},\ldots,\bar{\nu}_t^{h(\bar{\nu}_t)})$$ with $\bar{\nu}_1 > \cdots > \bar{\nu}_t,$
where $\bar{\nu}_i=\frac{\bar{r}_i}{\bar{s}_i}\in \Q .$ The stabilizer of $\nu$ with respect to the action of $S_d$ is
then given by
$$S_d(\nu):= S_{h(\bar{\nu}_1)} \times \ldots \times S_{h(\bar{\nu}_t)}.$$
Hence, $S_l(\nu):=S_l\cap S_d(\nu)$ is the stabilizer of $\nu$ with respect to the action of $S_l.$ If we define for $i=1,\ldots,t$ the number
$m_i\in \N$ as 
$$m_i:=\#\{j;\nu(j)=\bar{\nu}_i,\; j=1,\ldots,l\}$$
then we have $h(\bar{\nu}_i)= \bar{s}_i{\cdot} m_i$ for $ i=1,\ldots,t.$ We obtain
$$S_l(\nu)=S_{m_1} \times \cdots \times S_{m_t} \subset S_l.$$ 
Now we can certainly find a partition $(i_1,\ldots,i_r)$ of $l$ and an element $x\in S_l$ such that if we set
$k_j:=i_1+\cdots +i_j,$ we have for all $j=1,\ldots,r$
$$d_j=\sum_{k=k_{j-1}+1}^{k_j}s_{x^{-1}(k)}$$ 
and
\begin{eqnarray*}
\nu(W^\bullet,j) &  = & ({\nu({x^{-1}(k_{j-1}+1)}})^{s_{x^{-1}(k_{j-1}+1)}},\ldots,{\nu({x^{-1}(k_j)}})^{s_{x^{-1}(k_j)}} )\\
& = & (x\cdot\nu_{d_1+\cdots +d_{j-1}+1},\ldots, x\cdot\nu_{d_1+\cdots +d_j})  \in (\Q^{d_j})_+ \;.
\end{eqnarray*}
Here we put $k_0=0.$
Denote by $S_{i_1}\times \dots \times S_{i_r}$ the corresponding parabolic subgroup of $S_l.$ 
In the following we always identify double cosets with their corresponding Kostant representatives.
Recall that a Kostant representative $x\in S_l$ is the uniquely determined element in its double coset 
$$[x]\in  S_{i_1}\times \dots \times S_{i_r}  \backslash S_l/S_l(\nu)$$ such that
$x \mbox{ is strictly increasing on the intervals }$ 
$$[1,m_1],[m_1+1,m_1+m_2],\ldots,[m_1+\ldots + m_{t-1}+1,l]$$
and 
$x^{-1} \mbox{ is strictly increasing on the intervals }$ 
$$[1,i_1],[i_1+1,i_1+i_2],\ldots,[i_1+\ldots + i_{r-1}+1,l].$$

Now we fix as in \cite{R1} a system  $\underline{g}=(g_i)_{i=1,\ldots,r}$ of non-zero functions $g_i:\Q \rightarrow \Z_{\geq 0}$
such that $$\sum_{i=1}^r  g_i =g \;\mbox{ and }\; \sum_{y\in \Q} g_i(y)=d_i.$$
We associate to each $g_i$ its Hodge vector $\mu(g_i)\in (\Q^{d_i})_+$ by setting
$$\mu(g_i)=(\bar{\mu}_1^{g_i(\bar{\mu}_1)},\ldots,\bar{\mu}_s^{g_i(\bar{\mu}_s)})$$
(If some exponent vanishes then we omit it from the tuple).
Let $S_d(\mu)$ be the stabilizer of $\mu$  with respect to the action of $S_d$ on $\Q^d.$
It is given by
$$S_d(\mu)= S_{g(\bar{\mu}_1)} \times \ldots \times S_{g(\bar{\mu}_s)}.$$
We associate to the system of functions $\underline{g}$ a 
double coset  
$$[w]\in S_{d_1}\times \dots \times S_{d_r} \backslash  S_d/S_d(\mu),$$ 
with Kostant representative $w\in S_d$
such that for $i=1,\ldots,r$ 
$$\mu(g_i)=(w\cdot\mu_{d_1+ \cdots +d_{i-1}+1},\ldots, w\cdot\mu_{d_1+ \cdots +d_{i}}) .$$
We define the subset $\Fl(W^\bullet; \underline{g})$ of $\Fl(K)$ by
$$\Fl(W^\bullet; \underline{g}):=\{\F \in \Fl(K); \F|gr^i(W^\bullet) \in \Fl(gr^i(W^\bullet),\mu(g_i))^{ss}, i=1,\ldots,r \}.$$
The next lemma tells us when this set in non-empty. A necessary condition is of course the non-emptiness of 
$\Fl(gr^i(W^\bullet),\mu(g_i))^{ss}\;$ for all $i.$ It turns out that this condition is sufficient as well.

\begin{Lemma} We have
$\Fl(W^\bullet; \underline{g}) \neq \emptyset \mbox{ if and only if } \Fl(gr^i(W^\bullet),\mu(g_i))^{ss} \neq \emptyset$ for all $i=1,\ldots, r.$
\end{Lemma}

\proof  
Since the category of isocrystals over $k$ is semisimple (or by the definition of $W^\bullet$) we may choose a splitting $V=\oplus_{i=1}^r \tilde{W}_i$ of $W^\bullet.$
Pick for each $i\in \{1,\ldots, r\}$ a semistable filtration $\F_i \in \Fl(gr^i(W^\bullet), \mu(g_i))^{ss}=
\Fl(\tilde{W}_i,\mu(g_i))^{ss}.$  We get by building the sum $\oplus_{i=1}^r \F_i$ a filtration $\F \in \Fl(K)$, which lies by construction 
in $\Fl(W^\bullet; \underline{g}) .$ \qed

\noindent Thus we can make use of Lemma 3 to conclude when $\Fl(W^\bullet; \underline{g})$ is non-empty. 

\begin{Proposition}
We have $\Fl(W^\bullet; \underline{g}) \neq \emptyset $ if and only if $\mu(g_i) \geq \nu(W^\bullet,i)$
for all $ i=1,\ldots, r .$
\end{Proposition}

\begin{Remark} Suppose that the following strict inequalities are fulfilled for $i=1,\ldots,r-1$:
\begin{equation} 
\frac{|\mu(g_{i})-\nu(W^\bullet,i)|}{d_{i}} > \frac{|\mu(g_{i+1})-\nu(W^\bullet,i+1)|}{d_{i+1}}
\end{equation}
Then by construction each element in $\Fl(W^\bullet; \underline{g})$ has $W^\bullet$ as its HN-filtration.
\end{Remark}

Let $\C_K$ be the completion of an algebraic closure $\overline{K}$ of $K.$ This yields again an algebraically closed field.
We may define $\Fl(\C_K)^{ss}$
and $\Fl(W^\bullet,\underline{g})(\C_K)$ completely analogously as for $K.$
Then $\Fl(\C_K)^{ss}$ has the structure 
of a rigid-analytic variety (\cite{RZ}, Prop.1.36)  and is called the {\it period domain} with respect 
to $(V,\Phi,\mu)$. It is defined over $\Qps,$ where $s=lcm(s_1,\ldots,s_l).$ 
It is easily seen that the subsets $\Fl(W^\bullet,\underline{g})(\C_K)$  
are rigid-analytic varieties over $\Qps$ as well. 
Indeed, $\Fl(W^\bullet;\underline{g})(\C_K)$ is exactly the subset 
$\Fl(P(W^\bullet),\mu_P)$ of \cite{R2}. Here $P(W^\bullet)\subset J$ denotes the stabilizer of $W^\bullet,$
which is a parabolic subgroup of $J$ defined over $\Qp.$ 
Further we have an action of the $p$-adic group $P(W^\bullet)(\Qp)$ on  $\Fl(W^\bullet;\underline{g})(\C_K)$ 
in which the unipotent radical acts trivially.
We can then state:
\begin{Proposition}
The map $$\pi:\Fl(W^\bullet;\underline{g})(\C_K) \longrightarrow  
\prod_{i=1}^r \Fl(gr^i(W^\bullet),\mu(g_i))(\C_K)^{ss},$$ 
which sends a flag $\F$ to the restrictions $(\F|gr^i(W^\bullet))_{i=1,\ldots,r},$ is a $P(W^\bullet)(\Qp)$-equivariant vector bundle 
of rigid-analytic varieties of rank $l(w)$ over the base.
\end{Proposition}

We will show more generally the following statement. The proposition above is then simply an application of base change in the category
of rigid-analytic varieties.
Let $P_{GL(V)}(W^\bullet)$ be the stabilizer of the filtration $W^\bullet$ in $GL(V).$ Then we have
$P(W^\bullet)(\Qp)\subset P_{GL(V)}(W^\bullet)(K).$ Let $X$ be the locally closed subvariety 
of $\Fl$ (in the Zariski topology) consisting of flags $\F$ such that $\F|gr^i(W^\bullet)\in \Fl(gr^i(W^\bullet),\mu(g_i)) $ for $i=1,\ldots,r.$
Then we have:

\begin{Proposition}
The map $$\pi:X \longrightarrow  \prod_{i=1}^r \Fl(gr^i(W^\bullet),\mu(g_i)),$$ 
which sends a flag $\F$ to the restrictions $(\F|gr^i(W^\bullet))_{i=1,\ldots,r},$ is a $P_{GL(V)}(W^\bullet)$-equivariant vector bundle
 of varieties  of rank $l(w)$ over the base.
\end{Proposition}

\proof Put $P=P_{GL(V)}(W^\bullet).$ The morphism is clearly $P$-equivariant and surjective. Since the base is obviously smooth 
it is enough to show that the fibres are affine spaces of dimension $l(w).$
Let $Q=Q(\mu)$ be the parabolic subgroup which defines our flag variety $\Fl,$ i.e., $\Fl\cong GL(V)/Q.$
Then $X$ can be identified with the subvariety $PwQ/Q \;(\cite{O}\; 3.3).$ 
Let $P=M\cdot U$ be a Levi decomposition of $P.$ The right hand side of the morphism $\pi$ 
can be identified with 
$M/M\cap wQw^{-1}.$ The map $\pi$ is then the composition of the natural maps
$$ PwQ/Q \rightarrow P/P\cap wQw^{-1} \rightarrow M/M \cap wQw^{-1}.$$
The  fibres are thus isomorphic to 
$U/ U \cap wQw^{-1}.$ But this variety is well-known to be isomorphic to ${\mathbb A}^{l(w)}.$ \qed

Next we are interested in the union of those $\Fl(W^\bullet;\underline{g})$ where $W^\bullet$ varies over the set of filtrations
of $V$ by subisocrystals having the same Newton vectors as $W^\bullet.$ This union equals just
$$ \bigcup_{j\in J(\Qp)} \Fl(j\cdot W^\bullet;\underline{g}).$$
This set depends only on the Newton vectors of $W^\bullet$ and the Hodge vectors $\mu(g_i).$
For this reason we denote
this subset by $\Fl_\gamma,$ where the new parameter $\gamma$ 
is a triple $$\gamma=(P,[x],[w])$$ consisting of the partition $P=(i_1,\ldots,i_r)$ and
the double cosets 
$$[x]\in  S_{i_1}\times \dots \times S_{i_r}  \backslash S_l/S_l(\nu),$$  
$$[w]\in S_{d_1}\times  \dots \times S_{d_r}  \backslash  S_d/S_d(\mu)$$ 
defined above. For any such triple $\gamma=(P,[x_\gamma],[w_\gamma])$ we set for $j=1,\ldots,r$
\begin{itemize}
\item $d_j(\gamma):=\sum_{k=i_1+\cdots +i_{j-1}+1}^{i_1+\cdots +i_j}s_{x_\gamma^{-1}(k)} $  
\item $t_j(\gamma)=\sum_{k=1}^{j} d_\gamma(k)$  
\item $\nu(\gamma,j):= (x_\gamma\cdot\nu_{d_1+\cdots +d_{j-1}+1},\ldots, x_\gamma\cdot\nu_{d_1+\cdots +d_j})\in (\Q^{d_{j}(\gamma)})_+$
\item $\mu(\gamma,j):= (w_\gamma\cdot\mu_{d_1+ \cdots +d_{j-1}+1},\ldots, w_\gamma\cdot\mu_{d_1+ \cdots +d_{j}}) \in (\Q^{d_{j}(\gamma)})_+$
\item $\lambda(\gamma,j):=\mu(\gamma,j)-\nu(\gamma,j) \in \Q^{d_j(\gamma)}.$

\end{itemize}
Proposition 5 and Remark 6 lead to the following definition.
We define 
\begin{eqnarray*}
\Gamma:=\Gamma(V,\Phi,\mu) &:=&\Big\{\gamma=(P,[x_\gamma],[w_\gamma]); P=(i_1,\ldots,i_r) \mbox{ is a partition of } l ,\\
\\
& & \;\;[x_\gamma]\in  S_{i_1}\times \dots \times S_{i_r}  \backslash S_l/S_l(\nu),\\
\\
& & \;\; [w_\gamma]\in S_{d_1(\gamma)}\times \dots \times S_{d_r(\gamma)}  \backslash  S_d/S_d(\mu),\\
\\
& &\;\; \mbox{ such that } \lambda(\gamma,i)\geq 0  \mbox { in } \Q^{d_\gamma(i)} \mbox{ and }  \\
& & \\
& & \frac{|\lambda(\gamma,i)|}{d_{i}(\gamma)} > \frac{|\lambda(\gamma,i+1)|}{d_{i+1}(\gamma)}\; \forall i=1,\ldots,r-1 \Big\}.
\end{eqnarray*}
We obtain a disjoint union (a priori as sets)
$$ \Fl(K)=\bigcup^\cdot_{\gamma \in \Gamma} \Fl_\gamma.$$ 
and call the subsets $\Fl_\gamma$ the HN-strata with respect to $(V,\Phi,\mu).$
By the discussion above  $\Fl_\gamma$ is nothing else but
\begin{eqnarray*}
\Fl_\gamma=\Big\{ \F \in \Fl(K) & ;& |HN(\F)|=r, \dim gr^i(HN(\F))=d_i(\gamma), \\
& & \\
& & \F|gr^i(HN(\F)) \in \Fl(gr^i(HN(\F), \mu(\gamma,i)))^{ss}, \\
& & \\
& &  \nu(gr^i(HN(\F)))= \nu(\gamma,i), \; i=1,\ldots,r \Big\} .
\end{eqnarray*}
Similarly we may define $\Fl_\gamma(\C_K).$ It follows from Theorem 3' \cite{FR} that the resulting stratification of $\Fl(\C_K)$ has 
the same index set as  $\Fl(K).$ 
We want to stress that the sets $\Fl_\gamma(\C_K)$ are in general not rigid-analytic varieties as the following example demonstrates.
\begin{Example}
Let $d=2, \nu=(0,0) \mbox{ and } \mu$ arbitrary, i.e. we have $\Fl=\P^1$ and we are in the trivial isocrystal situation . 
Then we obtain as stratification 
$\P^1(\C_K)=(\P^1(\C_K)\setminus \P^1(K)) \; \bigcup\limits^\cdot \; \P^1(K).$ But $\P^1(K)$ is not a rigid-analytic variety.
\end{Example}

Since the sets $\Fl_\gamma(\C_K)$ are not well-behaved geometric objects, we may make use of Huber adic spaces \cite{H}.
Let $\Fl(j\cdot W^\bullet;\underline{g})^{ad}\subset \Fl^{ad}$ be the adic space corresponding to the rigid-analytic variety 
$\Fl(j\cdot W^\bullet;\underline{g}).$
Set
$$\Fl_\gamma^{ad}:= \bigcup_{j\in J(\Qp)} \Fl(j\cdot W^\bullet;\underline{g})^{ad},$$
where $\gamma$ corresponds to the data $W^\bullet,\underline{g}$ as above.
By definition $\Fl(j\cdot W^\bullet;\underline{g})^{ad}$ is a so-called prepseudo-adic subspace of $\Fl^{ad}$ (\cite{H} 1.10.1), i.e. a subset
of the adic space $\Fl^{ad}.$ By the specialization theorem for HN-polygons (\cite{R2} Theorem 3) we see that it is locally closed in $\Fl^{ad}.$ 
Moreover one sees easily that it is in fact a so-called pseudo-adic space which means that it is locally pro-constructible in the 
adic topology and convex with respect to the specializing order of points (\cite{H} 1.10.3). We obtain the first part of Theorem 2.

\begin{Theorem}
The pseudo-adic spaces $\Fl^{ad}_\gamma$ induce a stratification $\Fl^{ad}=\bigcup\limits^\cdot_{\gamma \in \Gamma}\Fl^{ad}_\gamma.$ 
\end{Theorem}

Put $P(\gamma)=P(W^\bullet)(\Qp).$ The second part of Theorem 2 is formulated in the  proposition below, which is an immediate consequence 
of Proposition 7.

\begin{Proposition}
The natural map 
$$\pi:\Fl_\gamma^{ad} \longrightarrow  J(\Qp)\times^{P(\gamma)}(\prod_{i=1}^r \Fl(gr^i(W^\bullet),\mu(g_i))^{ss})^{ad}$$ 
which is induced by the map of Proposition 7 is a $J(\Qp)$-equivariant vector bundle of pseudo-adic spaces 
of rank $l(w_{\gamma})$ over the base. 
\end{Proposition}

\begin{Remark}
The proposition above is the $p$-adic version of Proposition 2.6 \cite{R4}, which treats the finite field situation.
\end{Remark}

Now we want to give some examples.

\begin{Example} Consider the case of a trivial isocrystal, i.e., $V=K^d, \Phi=id_V \otimes \sigma.$ Then we have 
$\nu=(0,\ldots,0),$  $l=d$ and  $S_l(\nu)=S_l.$ The condition $\lambda(\gamma,\cdot)\geq 0$  is then automatically fulfilled.
Hence we get for any Hodge vector $\mu \in (\Q^d)_+:$
\begin{eqnarray*}
\Gamma= \Big\{\gamma=(P,[1],[w_\gamma])& ;&  P=(i_1,\ldots,i_r) \mbox{ a partition of $d$},\\
& & \\
& & [w_\gamma]\in S_{d_1(\gamma)}\times \dots \times S_{d_r(\gamma)} \backslash  S_d/S_d(\mu),\\
& & \\
& & \mbox{ such that for all } i=1,\ldots,r-1\\ 
& & \\
& & \sum_{j=t_{i-1}(\gamma)+1}^{t_i(\gamma)} \frac{w\cdot\mu_j}{d_i(\gamma)}  >  
\sum_{j=t_i(\gamma)+1}^{t_{i+1}(\gamma)}\frac{w\cdot\mu_j}{d_{i+1}(\gamma)}\Big\}.
\end{eqnarray*}
This is exactly the parametrisation of \cite{R1} p. 171.
\end{Example}

\begin{Example} 
Consider more generally the case of a basic isocrystal, i.e.,  
$\nu(1)=\nu(2)=\cdots =\nu(r).$ Again we have  $S_l(\nu)=S_l$ and the validity of the condition $\lambda(\gamma,\cdot)\geq 0.$  
Hence we get for any Hodge vector $\mu \in (\Q^d)_+:$
\begin{eqnarray*}
\Gamma= \Big\{\gamma=(P,[1],[w_\gamma])& ;&  P=(i_1,\ldots,i_r) \mbox{ a partition of $l$},\\
& & \\
& & [w_\gamma]\in S_{d_1(\gamma)}\times \dots \times S_{d_r(\gamma)} \backslash  S_d/S_d(\mu),\\
& & \\
& & \mbox{ such that for all } i=1,\ldots,r-1\\ 
& & \\
& & \sum_{j=t_{i-1}(\gamma)+1}^{t_i(\gamma)} \frac{w\cdot\mu_j}{d_i(\gamma)}  >  
\sum_{j=t_i(\gamma)+1}^{t_{i+1}(\gamma)}\frac{w\cdot\mu_j}{d_{i+1}(\gamma)}\Big\}.
\end{eqnarray*}
\end{Example}

\begin{Example}
Let $d=3,$ $\mu=(\mu_1>\mu_2>\mu_3)$ and $\nu=(\frac{1}{2},\frac{1}{2},0)$. Thus $\Fl$ is the complete flag variety
of $K^3$ and we have
\begin{itemize}
\item $l=2$,
\item $J(\Qp)=D^\ast\times \Qp^\ast$, where $D$ is a quaternion algebra over $\Qp$,
\item $S_\nu =S_2, S_2(\nu)=\{1\},$
\item $S_\mu=S_3,S_3(\mu)=\{1\}.$
\end{itemize}
The strata are indexed by the following elements of $\Gamma :$
$$\gamma_1=(P=(2),[1],[1]), \gamma_2=(P=(1,1),[1],[1]), \gamma_3=(P=(1,1),[(12)],[1])$$
Further we have $$\gamma_4=(P=(1,1),[1],[(23)])  \in  \Gamma   \Leftrightarrow  \mu_1 + \mu_3 > 2\mu_2 +1 .$$
\end{Example}

\vspace{1cm}
\begin{Remark} If we attach to $\gamma$ the polygon $P_\gamma$ in the euclidian plane with breaking points 
$$(0,0), (t_1(\gamma),|\lambda(\gamma,1)|), (t_2(\gamma),|\lambda(\gamma,1)| + |\lambda(\gamma,2)|) , \ldots $$
then condition $(3)$ of Remark 6 says that this polygon is convex. It is called the HN-polygon of $\gamma$.
Put $q_i:=\frac{|\lambda(\gamma,i)|}{d_i(\gamma)}.$ Then $P_\gamma$ has the shape:
\end{Remark}

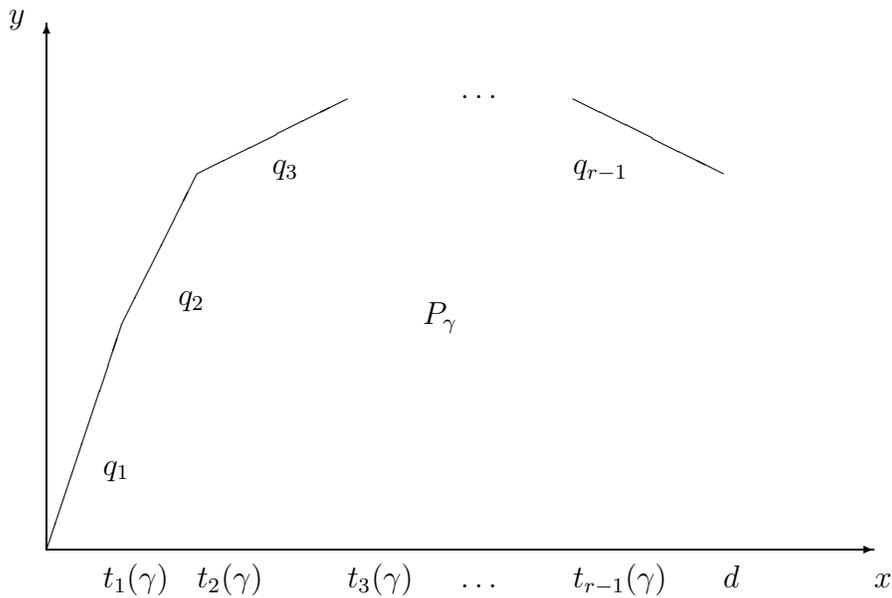
\setlength{\unitlength}{1cm}
\begin{picture}(12,8)
\put(0,0){\vector(0,1){7}}
\put(0,0){\vector(1,0){11}}

\put(11,-0.5){$x$}
\put(-0.5,7){$y$}
\put(5,3){$P_\gamma$}
\put(0.75,-0.5){$t_1(\gamma)$}
\put(2,-0.5){$t_2(\gamma)$}
\put(4,-0.5){$t_3(\gamma)$}
\put(5.5,-0.5){$\ldots$}
\put(7,-0.5){$t_{r-1}(\gamma)$}
\put(9,-0.5){$d$}
\put(0,0){\line(1,3){1}}
\put(1,3){\line(1,2){1}}
\put(2,5){\line(2,1){2}}
\put(7,6){\line(2,-1){2}}
\put(5.5,6){$\ldots$}
\put(0.75,1){$q_1$}
\put(1.75,3.25){$q_2$}
\put(3,5){$q_3$}
\put(7,5){$q_{r-1}$}

\end{picture}

\vspace{1.5cm}

\section{Theorem 1}
Now we turn to the question of Fontaine and Rapoport. 
For every $\gamma \in\Gamma$ we set
$$\lambda_\gamma=\Big(\Big(\frac{|\lambda(\gamma,1)|}{d_1(\gamma)}\Big)^{d_1(\gamma)},\Big(\frac{|\lambda(\gamma,2)|}{d_2(\gamma)}\Big)^{d_2(\gamma)},\ldots,\Big(\frac{|\lambda(\gamma,r)|}{d_r(\gamma)}\Big)^{d_r(\gamma)}\Big)\in (\Q^d)_+.$$ 
It is the HN-vector of any $x$ lying in $\Fl_\gamma.$ It coincides with the tuple $\lambda(x)$ of \cite{FR} Remark 3c.  
Thus the index set of all $\lambda(x)$ which appear in loc.cit. is the set 
$\Lambda:=\Lambda(V,\Phi,\mu):= \{\lambda_\gamma; \gamma \in \Gamma\}.$
We obtain a surjective map
$$\Gamma \rightarrow \Lambda,$$
which is in general not injective as the following example shows.

\begin{Example} Let $d=5$, $\nu=(0,0,0,0,0), \mu=(\mu_1>\mu_2>\mu_3>\mu_4>\mu_5) \in (\Q^d)_+$ with 
$\mu_1+\mu_4=\mu_2+\mu_3=-\frac{\mu_5}{2}.$ Then the following two elements of $\Gamma$ 
$$\gamma_1=((2,3),[1],(2,3,4)) \mbox{ and }  \gamma_2=((2,3),[1],(1,3,2))$$
have the same HN-vector $\lambda=((\frac{\mu_1+\mu_4}{2})^2,(\frac{\mu_2+\mu_3+\mu_5}{3})^3).$
\end{Example}

\noindent We thus see that the stratification $\Fl(K)=\bigcup\limits^\cdot_{\gamma\in \Gamma}  \Fl_\gamma$ is finer than the stratification induced 
by the HN-vectors used in \cite{FR}. Using the surjective map $\Gamma \longrightarrow \Lambda$ we can describe the elements
of $\Q^d$ for which the corresponding strata with respect to $\Lambda$ are non-empty. This is the content of Theorem 1, which follows 
now easily.

\bigskip
\begin{flushleft}
Universit\"at Leipzig \\ 
Fakult\"at f\"ur Mathematik und Informatik \\ 
Mathematisches Institut \\ 
Augustusplatz 10/11 \\ 
D-04109 Leipzig \\ 
Germany 
\end{flushleft}

\end{document}